\documentclass[conference]{IEEEtran}
\usepackage{cite}
\usepackage{amsmath,amssymb,amsfonts}
\usepackage{algorithmic}
\usepackage{graphicx}
\usepackage{textcomp}
\usepackage{xcolor}
\usepackage{pifont}
\usepackage{hyperref}
\usepackage{cleveref}
\usepackage{xcolor}
\usepackage{mathtools, cuted}
\usepackage{flushend}
\def\BibTeX{{\rm B\kern-.05em{\sc i\kern-.025em b}\kern-.08em
    T\kern-.1667em\lower.7ex\hbox{E}\kern-.125emX}}
\begin{document}

\title{A New Dynamic Model of a Two-Wheeled Two-Flexible-Beam Inverted Pendulum Robot\\
}

\author{Amin Mehrvarz$^{1}$*, Mohammad Javad Khodaei$^{1}$*, William Clark$^1$, and Nader Jalili$^{1,2}$ \\
$^1$Department of Mechanical and Industrial Engineering, Northeastern University, Boston, MA, U.S. \\
$^2$Professor and Head, Department of Mechanical Engineering, University of Alabama, Tuscaloosa, AL, U.S. \\
* Co-first author}

\maketitle

\begin{abstract}
Two-wheeled inverted pendulum robots are designed for self-balancing and they have remarkable advantages. In this paper, a new configuration and consequently dynamic model of one specific robot is presented and its dynamic behavior is analyzed. In this model, two cantilever beams are on the two-wheeled base and they are excited by voltages to the attached piezoelectric actuators. The mathematical model of this system is obtained using the extended Hamilton's Principle. The results show that the governing equations of motion are highly nonlinear and contain several coupled terms. These complex equations are solved numerically and the natural frequencies of the system are extracted. The system is then simulated in both lateral and horizontal plan movements. This paper proposes a new model of two-wheeled two-flexible-beam inverted pendulum Robot and investigates its complex dynamic; however, the derived equations will be validated experimentally and a suitable control strategy will be applied to the system to make it fully automated and more applicable in the future works.
\end{abstract}

\begin{IEEEkeywords}
Two-wheeled robot, flexible system, Euler-Bernoulli beam, inverted pendulum.
\end{IEEEkeywords}

\section{Introduction}
Inverted pendulums are traditional dynamic problems. If an inverted pendulum is used in a moving cart, new type of interesting problems will appear. One of these problems is two-wheeled inverted pendulum systems. Because of their small size, great performance in quick driving, and their stability with controller, scientists and engineers are interested in them\textcolor{blue}{\cite{jeong2018development}}. Based on this interest, each year new models are introduced and new robots are made. Self-transportation systems such as hoverboards and small two-wheeled robots are the most important patents based on the moving inverted pendulums\textcolor{blue}{\cite{canete2015modeling}}.

The idea of using self-transportation systems comes from a push in the transportation industry to develop transportation systems that contribute less to pollution and cause less damage to the environment overall. One approach to this has been a shift to Personal Electric Vehicles (PEVs), which are powered by electricity rather than combustion. PEVs provide many benefits to both consumers and society, including lower costs than automobiles, shorter trip times for short distances, cleaner transportation, and mobility for the disabled\textcolor{blue}{\cite{ulrich2005estimating}}. One popular type of PEV that has emerged in recent years is the "Stand-on Scooter". In 2005, Ulrich analyzed existing stand-on scooter technology and estimated that the light design of these PEVs combined with their modest range and speed would be "highly feasible technically, and with substantial consumer demand could be feasible economically"\textcolor{blue}{\cite{ulrich2005estimating}}. One type of stand-on scooter analyzed by Ulrich was the Segway, which is a type of PEV that marketed itself based on its inverted pendulum balancing mechanism and its agility.

Inverted-pendulum transporter is a type of self-balancing system that allows for an operator to control it without the need for a throttle. Instead, the device applies a lateral movement to the system based on an angle applied by the operator, who acts as an inverted pendulum, in order to keep the inverted pendulum balanced and stable in the upright position. A popular inverted pendulum PEV is the "Hoverboard", which consists of two motorized wheels connected to two independent articulating pads. The operator controls the speed of travel by leaning forward and backward and controls the direction of travel by twisting the articulating pads with their feet. This allows for an inexpensive transportation system that is compact, as it does not require any large motor or steering apparatus, and agile, as it can be controlled easily and move quickly. Considerable work has already been done in developing a transportation device that uses this sort of self-balancing mechanism. Grasser et al.  developed a scaled-down prototype of a self-balancing pendulum, but not a full-scale version that could be ridden by a person\textcolor{blue}{\cite{grasser2002joe}}. Tsai et al.\textcolor{blue}{\cite{tsai2010adaptive}} developed a self-balancing PEV but used handlebars for guidance rather than the articulating pads like the hoverboard.

Moving inverted pendulum systems are also used in robots for different objects. For example, Solis et al. used the robot for educational purposes\textcolor{blue}{\cite{solis2009development}} or Double Robotics Inc made a robot for telecommuter\textcolor{blue}{\cite{double2019}}. Two steps are needed to design a proper controller for these robots. Finding a perfect model for the two-wheeled inverted pendulum is the first step and designing the best controller is the second step. It should be mentioned that the self-balancing systems like these robots create difficult control problems, as they are inherently unstable and subject to unpredictable external forces from their environment. To solve these problems, a robust dynamic model must be developed, in order to fully understand the physical properties of the system, and also a better comprehending of how to control it to maintain stable. These robots also need one or two motors for moving which their motors’ models should be considered in the robot’s dynamic\textcolor{blue}{\cite{frankovsky2017modeling}}.

Decreasing the number of wheels which are unnecessary in most of the time and only used for the system's balancing was one of the first ideas of the two-wheeled robots\textcolor{blue}{\cite{kim2005dynamic}}. For this purpose, Kim et al. developed a mathematical model for a self-balancing two-wheeled robot that was capable of changing direction. This robot acted as a rigid single-pendulum, allowing the model to assume a lumped-parameter system and they designed a linear controller for their robot \textcolor{blue}{\cite{kim2005dynamic}}. Other researchers tried to develop\textcolor{blue}{\cite{jeong2018development}} or modify the previous models\textcolor{blue}{\cite{kim2015dynamic}}. They also tried to improve the controller performance during the robots' operation. In a research, Zafar et al. discussed a derived mathematical model for a self-balancing inverted-pendulum robot and implement the operational space controller\textcolor{blue}{\cite{zafar2016whole}}. In other researches, new controllers for the two-wheeled inverted pendulum systems were designed. These controllers can work with time-varying parametric uncertainties\textcolor{blue}{\cite{boukens2017design}}, strong nonlinear behaviors due to abrupt external disturbances\textcolor{blue}{\cite{kim2017nonlinear}}, and initial errors, pulse disturbance and random noises\textcolor{blue}{\cite{yue2018efficient}}.

A new type of two-wheeled inverted pendulum robot is also designed for specific applications. This pendulum was not rigid and it behaves as a flexible system. Partial Differential Equations (PDEs) are the most standard way of mathematically representing continuous systems and has been used for vibrational analysis and control for different models and applications\textcolor{blue}{\cite{mehrvarz2019vibration, mehrvarz2018vibration, karagiannis2019exponential,marzban2016effect}}. Researchers work on the flexible structure such as beams and bars with the base motion for several years but they did not consider their systems as a moving robots. They used flexible structures with moving base for different applications, such as micro gyroscopes and piezoresponse force microscopy and etc.\textcolor{blue}{\cite{salehi2009vibration, ansari2009coupled, khodaei2018theoretical}}. The base motion creates different accelerations and complex nonlinear PDE equations of motion for the flexible system. The first idea of using flexible structure as an inverted pendulum in robot was created by Nguyen et al. They studied a linearized mathematical model and controller for a single flexible inverted pendulum, this model only accounted for lateral movement in one direction rather than in the two directions and it has four wheels\textcolor{blue}{\cite{nguyen2016designing}}. In another robot Mehrvarz et al. modeled a two-wheeled flexible inverted pendulum which can move in one direction\textcolor{blue}{\cite{mehrvarz2018modeling}}. They also designed a MPC controller for their robot and showed that the percision of their modeling\textcolor{blue}{\cite{clark2019control}}. Their robot cannot move in the plane and this was their main problem. 

Here, to solve the problem of two-wheeled flexible beam inverted pendulum which is used in\textcolor{blue}{\cite{mehrvarz2018modeling}}, two-wheeled two-flexible beam inverted pendulum model is addressed. This robot is designed to move in-plane and has not the previous problem. Because of its complexity and non-linearity, the dynamic model and the vibrations of the system are just considered and analyzed in this paper and controller is not designed for this system right now. The proposed model analyzes the pure dynamics of the robot and the vibrations of the beams at the same time, which is a novel approach. The main goal of this paper is to investigate and simulate the dynamic model of a piezoelectrically actuated cantilever beams on the two-wheeled robot.

The remainder of this paper is organized as follow. The dynamic equations of the system are derived in \textcolor{blue}{{Section }\ref{Mathematical modeling}} and a brief summary of their solution is presented in \textcolor{blue}{{Section }\ref{numericalsimulation}} In \textcolor{blue}{{Section }\ref{simulationresults}} the simulation results are discussed and a conclusion is given in \textcolor{blue}{{Section }\ref{conclusion}}.  

\section{Mathematical modeling }\label{Mathematical modeling}
In this section, the governing equations of motion are derived using the extended Hamilton's principle. In order to apply the extended Hamilton's principle to the system, the positions and translational and rotational velocities of these elements are needed to be defined mathematically. As shown in \textcolor{blue}{{Fig.}\ref{fig1}}, the system has two beams and two piezoelectric actuators to excite the beams. These beams are mounted on two independent bases, which are attached to two wheels and DC motors. In this model, as seen in \textcolor{blue}{{Fig.}\ref{fig1}}, two flexible cantilever beams act as flexible inverted pendulums fixed to two articulating bases. The mathematical model for the system represents the response of the system to small disturbances from piezoelectric actuators mounted to the base of each pendulum. These actuators cause deformation and bending in the continuous pendulums when voltage is applied to them. As seen in \textcolor{blue}{{Fig.}\ref{fig2}}, the position of the right (beam$_{1}$) and left (beam$_2$) beams in the $XYZ$ frame are as:
\begin{figure}
\centering
\includegraphics[width=80mm]{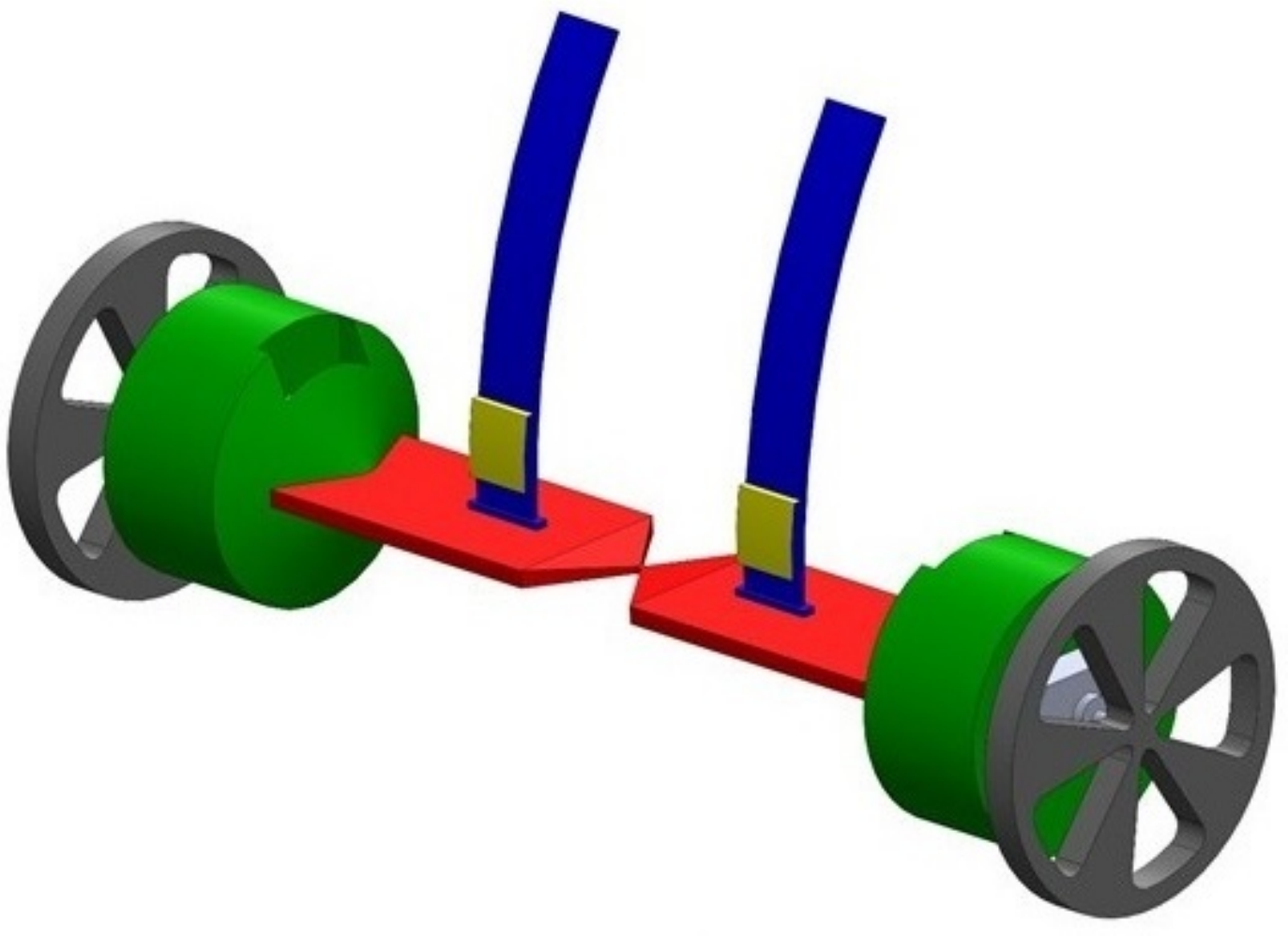}
\caption{ The proposed dynamic model of the two-wheeled two-flexible-beam inverted pendulum robot.}
\label{fig1}
\end{figure}

\begin{align}
\label{eq1}
\left[ \begin{matrix}
   {{X}_{beam{}_{1}}}  \\
   {{Y}_{beam{}_{1}}}  \\
   {{Z}_{beam{}_{1}}}  \\
\end{matrix} \right]=\left[ \begin{matrix}
   {{X}_{c}}+a\sin \varphi +\cos \varphi \left( \cos {{\theta }_{1}}{{w}_{1}}+{{x}_{1}}\sin {{\theta }_{1}} \right)  \\
   {{Y}_{c}}-a\cos \varphi +\sin \varphi \left( \cos {{\theta }_{1}}{{w}_{1}}+{{x}_{1}}\sin {{\theta }_{1}} \right)  \\
   {{r}_{w}}-\sin {{\theta }_{1}}{{w}_{1}}+{{x}_{1}}\cos {{\theta }_{1}}  \\
\end{matrix} \right]
\end{align}

\begin{align}
\label{eq2}
\left[ \begin{matrix}
   {{X}_{beam{}_{2}}}  \\
   {{Y}_{beam{}_{2}}}  \\
   {{Z}_{beam{}_{2}}}  \\
\end{matrix} \right]=\left[ \begin{matrix}
   {{X}_{c}}-a\sin \varphi +\cos \varphi \left( \cos {{\theta }_{2}}{{w}_{2}}+{{x}_{2}}\sin {{\theta }_{2}} \right)  \\
   {{Y}_{c}}+a\cos \varphi +\sin \varphi \left( \cos {{\theta }_{2}}{{w}_{2}}+{{x}_{2}}\sin {{\theta }_{2}} \right)  \\
   {{r}_{w}}-\sin {{\theta }_{2}}{{w}_{2}}+{{x}_{2}}\cos {{\theta }_{2}}  \\
\end{matrix} \right]
\end{align}
where $X_c$ and $Y_c$ are the positions of the center of gravity of the robot and $a$ denotes the distance between the center of the beams and the center of gravity in the $y$-direction. The robot can rotate around the $z$-axis and the bases have different rotations around the $y$-axis. These rotational angles are shown by $\varphi$, $\theta_1$ and $\theta_2$.  In \eqref{eq1} and \eqref{eq2}, the parameters $w_1$ and $w_2$ represent the bending deflections of the beams and $r_w$ is the radius of the wheels. Since the beams are supposed to be continuous, the position of each particle in the beams are given by $x_1$ and $x_2$. The position of the right (wheel$_1$ and base$_1$) and left (wheel$_2$ and base$_2$) wheels and bases can be obtained as: 

\begin{equation}
\label{eq3}
\left[ \begin{matrix}
   {{X}_{whee{{l}_{1}}}}  \\
   {{Y}_{whee{{l}_{1}}}}  \\
   {{Z}_{whee{{l}_{1}}}}  \\
\end{matrix} \right]=\left[ \begin{matrix}
   {{X}_{c}}+2a\sin \varphi   \\
   {{Y}_{c}}-2a\cos \varphi   \\
   {{r}_{w}}  \\
\end{matrix} \right]  
\end{equation}

\begin{equation}
\label{eq4}
\left[ \begin{matrix}
   {{X}_{whee{{l}_{2}}}}  \\
   {{Y}_{whee{{l}_{2}}}}  \\
   {{Z}_{whee{{l}_{2}}}}  \\
\end{matrix} \right]=\left[ \begin{matrix}
   {{X}_{c}}-2a\sin \varphi   \\
   {{Y}_{c}}+2a\cos \varphi   \\
   {{r}_{w}}  \\
\end{matrix} \right]    
\end{equation}

\begin{equation}
 \label{eq5}
\left[ \begin{matrix}
   {{X}_{bas{{e}_{1}}}}  \\
   {{Y}_{bas{{e}_{1}}}}  \\
   {{Z}_{bas{{e}_{1}}}}  \\
\end{matrix} \right]=\left[ \begin{matrix}
   {{X}_{c}}+a\sin \varphi   \\
   {{Y}_{c}}-a\cos \varphi   \\
   {{r}_{w}}  \\
\end{matrix} \right]   
\end{equation}

\begin{equation}
\label{eq6}
\left[ \begin{matrix}
   {{X}_{bas{{e}_{2}}}}  \\
   {{Y}_{bas{{e}_{2}}}}  \\
   {{Z}_{bas{{e}_{2}}}}  \\
\end{matrix} \right]=\left[ \begin{matrix}
   {{X}_{c}}-a\sin \varphi   \\
   {{Y}_{c}}+a\cos \varphi   \\
   {{r}_{w}}  \\
\end{matrix} \right]
\end{equation}

Consequently, the translational velocity of all the elements can be obtained by applying a time derivation operator to \eqref{eq1} through \eqref{eq6}.

It is assumed that the system move without slippage in $y$ and $x$ directions. Hence, velocity of the center of gravity of the system can have the following relationship,
\begin{equation}\label{eq7}
    \frac{{{{\dot{Y}}}_{c}}}{{{{\dot{X}}}_{c}}}=-\tan \varphi 
\end{equation}

It should be noted that the whole velocity of the beams can be calculated through integration of the each beam's particles velocity along the beam. Besides the translational motion, the elements of the robot have also some rotational movements. The rotational velocity of each element can be written as: 

\begin{equation}
    \label{eq8}
    \left[ \begin{matrix}
   {{\omega }_{{{x}_{_{whee{{l}_{1}}}}}}}  \\
   {{\omega }_{{{y}_{_{whee{{l}_{1}}}}}}}  \\
   {{\omega }_{{{z}_{_{whee{{l}_{1}}}}}}}  \\
\end{matrix} \right]=\left[ \begin{matrix}
   0  \\
   \frac{{{{\dot{X}}}_{c}}\cos \varphi +{{{\dot{Y}}}_{c}}\sin \varphi +2a\dot{\varphi }}{{{r}_{w}}}  \\
   {\dot{\varphi }}  \\
\end{matrix} \right]
\end{equation}

\begin{equation}
    \label{eq9}
    \left[ \begin{matrix}
   {{\omega }_{{{x}_{_{whee{{l}_{2}}}}}}}  \\
   {{\omega }_{{{y}_{_{whee{{l}_{2}}}}}}}  \\
   {{\omega }_{{{z}_{_{whee{{l}_{2}}}}}}}  \\
\end{matrix} \right]=\left[ \begin{matrix}
   0  \\
   \frac{{{{\dot{X}}}_{c}}\cos \varphi +{{{\dot{Y}}}_{c}}\sin \varphi -2a\dot{\varphi }}{{{r}_{w}}}  \\
   {\dot{\varphi }}  \\
\end{matrix} \right]
\end{equation}

\begin{equation}
   \label{eq10}
    \left[ \begin{matrix}
   {{\omega }_{{{x}_{_{bas{{e}_{1}}}}}}}  \\
   {{\omega }_{{{y}_{_{bas{{e}_{1}}}}}}}  \\
   {{\omega }_{{{z}_{_{bas{{e}_{1}}}}}}}  \\
\end{matrix} \right]=\left[ \begin{matrix}
   0  \\
   {{{\dot{\theta }}}_{1}}  \\
   {\dot{\varphi }}  \\
\end{matrix} \right] 
\end{equation}

\begin{equation}
    \label{eq11}
\left[ \begin{matrix}
   {{\omega }_{{{x}_{_{bas{{e}_{2}}}}}}}  \\
   {{\omega }_{{{y}_{_{bas{{e}_{2}}}}}}}  \\
   {{\omega }_{{{z}_{_{bas{{e}_{2}}}}}}}  \\
\end{matrix} \right]=\left[ \begin{matrix}
   0  \\
   {{{\dot{\theta }}}_{2}}  \\
   {\dot{\varphi }}  \\
\end{matrix} \right] 
\end{equation}

\begin{equation}
    \label{eq12}
    \left[ \begin{matrix}
   {{\omega }_{{{x}_{_{bea{{m}_{1}}}}}}}  \\
   {{\omega }_{{{y}_{_{bea{{m}_{1}}}}}}}  \\
   {{\omega }_{{{z}_{_{bea{{m}_{1}}}}}}}  \\
\end{matrix} \right]=\left[ \begin{matrix}
   \dot{\varphi }\cos {{\theta }_{1}}  \\
   {{{\dot{\theta }}}_{1}}+\frac{{{\partial }^{2}}{{w}_{1}}}{\partial x\partial t}  \\
   \dot{\varphi }\sin {{\theta }_{1}}  \\
\end{matrix} \right]
\end{equation}

\begin{equation}
    \label{eq13}
   \left[ \begin{matrix}
   {{\omega }_{{{x}_{_{bea{{m}_{2}}}}}}}  \\
   {{\omega }_{{{y}_{_{bea{{m}_{2}}}}}}}  \\
   {{\omega }_{{{z}_{_{bea{{m}_{2}}}}}}}  \\
\end{matrix} \right]=\left[ \begin{matrix}
   \dot{\varphi }\cos {{\theta }_{2}}  \\
   {{{\dot{\theta }}}_{2}}+\frac{{{\partial }^{2}}{{w}_{2}}}{\partial x\partial t}  \\
   \dot{\varphi }\sin {{\theta }_{2}}  \\
\end{matrix} \right]
\end{equation}

Since the system has 6 different parts, the kinetic energy of the whole system including the translational and rotational parts can be obtained as:

\begin{equation}\label{eq14}
T=\frac{1}{2}\sum\limits_{i=1}^{6}{\left( {{\rho }_{i}}{{A}_{i}}{{V}_{i}}^{2}+{{I}_{xi}}{{\omega}_{xi}}^{2}+{{I}_{yi}}{{\omega}_{yi}}^{2}+{{I}_{zi}}{{\omega }_{zi}}^{2} \right)}
\end{equation}
where ${V_i}^2={\dot{X_i}}^2+{\dot{Y_i}}^2+{\dot{Z_i}}^2$. $I_{xi}$, $I_{yi}$ and $I_{zi}$ are the mass moments of inertia of the $i$-th element and are assumed to be equal for each wheel, each base, and each beam. Also, the effect of rotary inertia terms of the beams are ignored as in\textcolor{blue}{\cite{bhadbhade2008novel}}. The combined $\rho A$ for the system can be calculated as:

\begin{equation}\label{eq15}
\rho A=\left\{ \begin{matrix}
   {{\rho }_{b}}{{A}_{b}}+{{\rho }_{p}}{{A}_{p}} & 0<x\le {{L}_{p}}  \\
   {{\rho }_{b}}{{A}_{b}} & {{L}_{p}}<x\le{L}  \\
\end{matrix} \right.
\end{equation}
where $L$ and $L_p$ are the beam and the piezoelectric lengths, $\rho_b$ and $\rho_p$ are the densities of the beam and piezoelectric actuators, respectively, and $A_b$ and $A_p$ are the cross-sectional areas. Also, the potential energy of the beams and piezoelectric actuators can be expressed as:

\begin{equation}\label{eq16}
\begin{split}
  & U=2\rho gAL+{{\int_{0}^{L}{\frac{1}{2}{{E}_{b}}{{I}_{b}}\left( \frac{{{\partial }^{2}}{{w}_{1}}}{\partial {{x}_{1}}^{2}} \right)}}^{2}}d{{x}_{1}}\\
  &+{{\int_{0}^{{{L}_{p}}}{\frac{1}{2}{{E}_{p}}{{I}_{p}}\left( \frac{{{\partial }^{2}}{{w}_{1}}}{\partial {{x}_{1}}^{2}}+{{z}_{p}}{{d}_{31}}\frac{{{v}_{1}}\left( t \right)}{{{t}_{p}}} \right)}}^{2}}d{{x}_{1}}\\
  &+{{\int_{0}^{L}{\frac{1}{2}{{E}_{b}}{{I}_{b}}\left( \frac{{{\partial }^{2}}{{w}_{2}}}{\partial {{x}_{2}}^{2}} \right)}}^{2}}d{{x}_{2}}\\
  &+{{\int_{0}^{{{L}_{p}}}{\frac{1}{2}{{E}_{p}}{{I}_{p}}\left( \frac{{{\partial }^{2}}{{w}_{2}}}{\partial {{x}_{2}}^{2}}+{{z}_{p}}{{d}_{31}}\frac{{{v}_{2}}\left( t \right)}{{{t}_{p}}} \right)}}^{2}}d{{x}_{2}} \\ 
  & +\int_{0}^{L}{\rho gA}\left( -\sin {{\theta }_{1}}{{w}_{1}}+{{x}_{1}}\cos {{\theta }_{1}} \right)d{{x}_{1}} \\
  &+\int_{0}^{L}{\rho gA}\left( -\sin {{\theta }_{2}}{{w}_{2}}+{{x}_{2}}\cos {{\theta }_{2}} \right)d{{x}_{2}} 
  \end{split}
\end{equation}
where $E_b$ and $E_p$ denote Young's modulus of elasticity of the beam and the piezoelectric, respectively, and  $I_b$ and $I_p$ are the mass moments of inertia of the beam and piezoelectric cross-section about $y$-axis, respectively. In Eq. \eqref{eq16}, $z_p$ is the neutral axis along the $z$-axis, $d_{31}$ denotes the piezoelectric constant of the actuator, $v_1 (t)$ and $v_2 (t)$ are voltages that are applied to the piezoelectric actuators and $t_p$ is the thickness of the piezoelectric actuators. 

The damping effects of the beams can be taken to account as virtual work terms as follows:  

\begin{equation}\label{eq17}
\delta {{W}^{nc}}=\int_{0}^{L}{{{C}_{1}}\frac{\partial {{w}_{1}}}{\partial t}\delta {{w}_{1}}}d{{x}_{1}}+\int_{0}^{L}{{{C}_{2}}\frac{\partial {{w}_{2}}}{\partial t}\delta {{w}_{2}}}d{{x}_{2}}
\end{equation}
where $C1$ and $C2$ are the viscous damping coefficients of the beams. As noted, the robot is assumed to have two DC motors, which produce torques $\tau_1$ and $\tau_2$. The work of these external torques makes additional virtual work terms as:

\begin{equation}\label{eq18}
\begin{split}
 &\delta {{W}^{ext}}=\frac{{{\tau }_{1}}+{{\tau }_{2}}}{{{r}_{w}}}\left(\right. \cos \varphi \delta {{X}_{c}}+\sin \varphi \delta Y+\left(\right. {{Y}_{c}}\cos \varphi\\
 &-{{X}_{c}}\sin \varphi  \left.\right)\delta \varphi  \left.\right)+{{F}_{s}}\left(\right. -\sin \varphi \delta {{X}_{c}}+\cos \varphi \delta {{Y}_{c}}-\\
 &\left(\right. {{X}_{c}}\cos \varphi +{{Y}_{c}}\sin \varphi  \left.\right)\delta \varphi  \left.\right) 
 \end{split}
\end{equation}

\begin{figure}
\centering
\includegraphics[width=80mm]{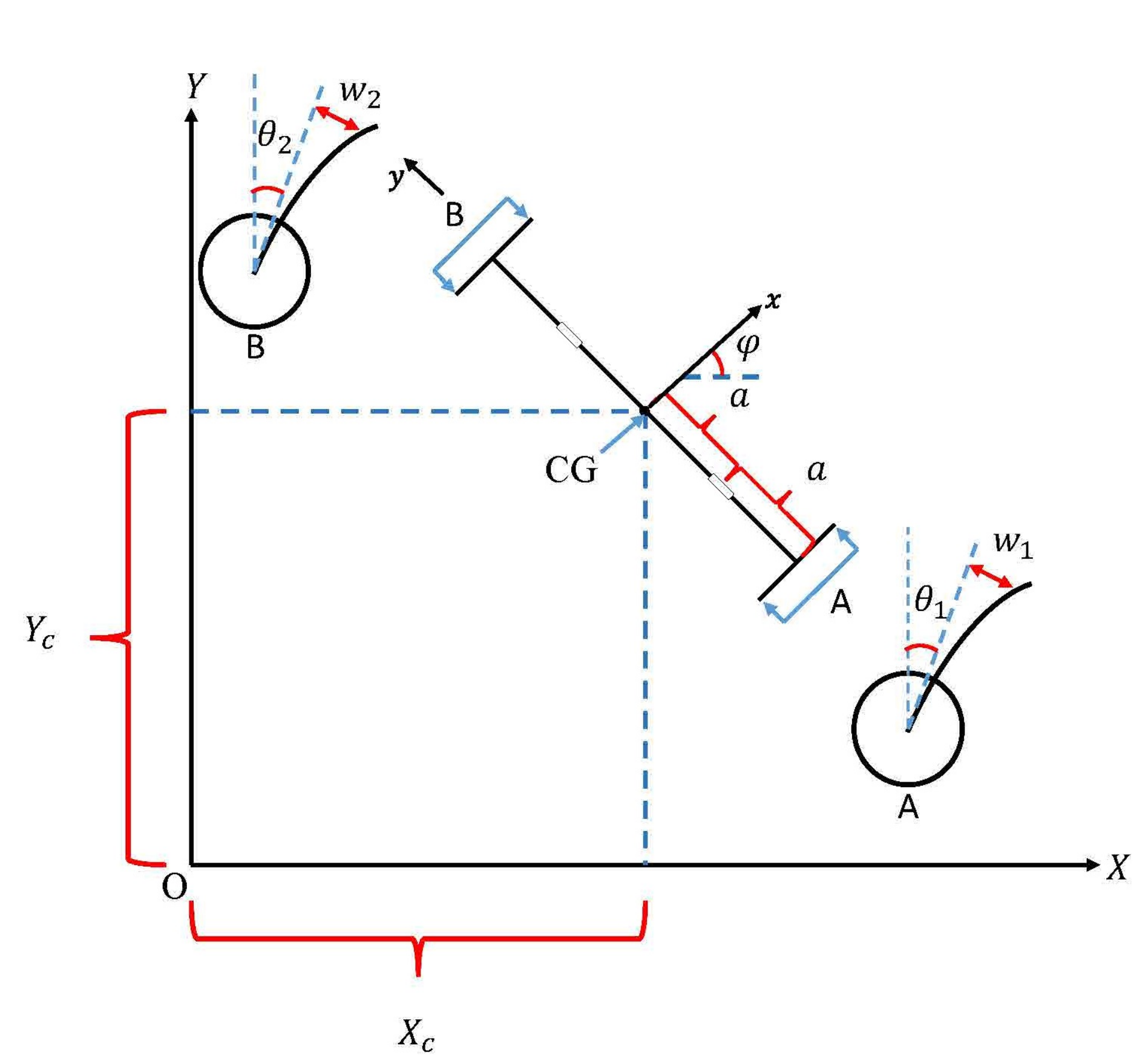}
\caption{ The robot kinematics.}
\label{fig2}
\end{figure}

\begin{figure}
\centering
\includegraphics[width=80mm]{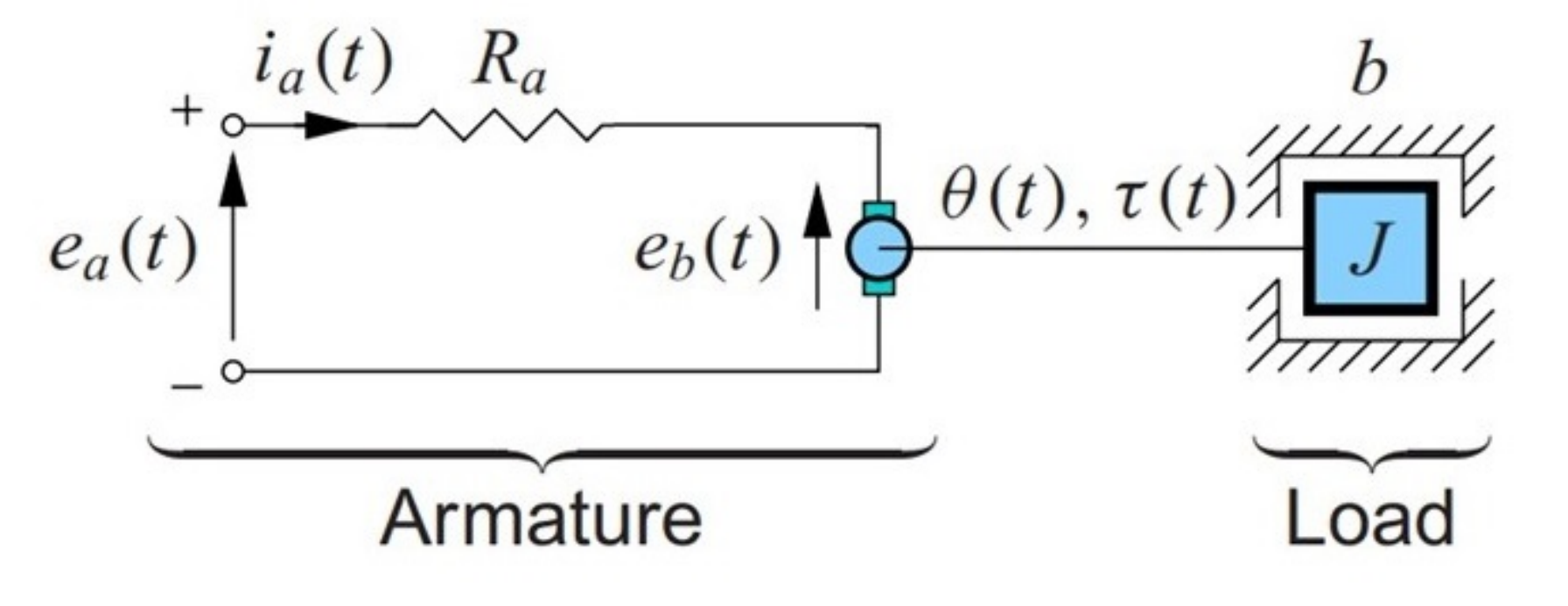}
\caption{The DC motor circuit diagram.}
\label{fig3}
\end{figure}

In Eq. \eqref{eq18}, the torques $\tau_1$ and $\tau_2$ are produced by the wheel-connected DC motors. The circuit diagram of the DC motors is shown in \textcolor{blue}{{Fig.}\ref{fig3}}. In this figure, the parameter $R_a$ denotes the armature resistance and the variables $V_a$ and $i_a$ are the applied voltage and the motor current draw, respectively. 
The equation of this system can be derived by applying the Kirchhoff's voltage Law to the circuit as:

\begin{equation}\label{eq19}
\begin{matrix}
   {{V}_{aj}}(t)={{R}_{a}}{{i}_{aj}}(t)+{{\tau }_{Bj}} & ,j=1,2  \\
\end{matrix}
\end{equation}
where $\tau_{Bi}$ represents the back electromotive force (emf) and is equal to:

\begin{equation}\label{eq20}
\begin{matrix}
   {{\tau }_{Bj}}={{K}_{B}}\frac{{{{\dot{X}}}_{whee{{l}_{j}}}}}{{{r}_{w}}} & ,j=1,2  \\
\end{matrix}
\end{equation}
with $K_B$ being the motor speed coefficient. Here,  torques $\tau_1$ and $\tau_2$ are given as:

\begin{equation}\label{eq21}
    \begin{matrix}
   {{\tau }_{j}}={{K}_{t}}{{i}_{aj}} & ,j=1,2  \\
\end{matrix}
\end{equation}
where $K_t$ is the motor torque constant and is provided by the manufacturer. Hence, the coupled equation can be obtained for the DC motor part by substituting Eq. \eqref{eq20} into Eq. \eqref{eq19}.

\begin{equation}\label{eq22}
\begin{matrix}
   {{V}_{aj}}(t)={{R}_{a}}{{i}_{aj}}(t)+{{K}_{B}}\frac{{{{\dot{X}}}_{whee{{l}_{j}}}}}{{{r}_{w}}} & ,j=1,2  \\
\end{matrix}
\end{equation}

The extended Hamilton's principle for this system can be written as: 

\begin{equation}\label{eq23}
\int_{0}^{t}{(\delta T-\delta U+\delta {{W}^{nc}}+\delta {{W}^{ext}})dt}=0
\end{equation}

After substituting Eqs. \eqref{eq14}-\eqref{eq18} to Eq. \eqref{eq23} and some manipulations and simplifications, the dynamic equations of the system can be obtained as:

\begin{strip}
\begin{align}
    \label{eq24}
\begin{split}
  & \left( \frac{{{\tau }_{1}}+{{\tau }_{2}}}{{{r}_{w}}} \right)\cos \varphi -\sin \varphi {{F}_{s}}+\frac{2{{I}_{{{y}_{wheel}}}}}{{{r}_{w}}^{2}}\left(\right. \sin 2\varphi \dot{X}\dot{\varphi }-\cos 2\varphi \dot{Y}\dot{\varphi } \left.\right)-2\left(\right. \frac{{{I}_{{{y}_{wheel}}}}}{{{r}_{w}}^{2}}{{\cos }^{2}}\varphi +{{m}_{w}}\\
  &+{{m}_{base}} \left.\right)\ddot{X}+\frac{{{I}_{{{y}_{wheel}}}}}{{{r}_{w}}^{2}}\sin 2\varphi \ddot{Y}-\int_{0}^{L}{\rho A}\left[\right. \ddot{X}+\left( a\cos \varphi -\cos {{\theta }_{1}}\sin \varphi {{w}_{1}}-{{x}_{1}}\sin {{\theta }_{1}}\sin \varphi  \right)\\
  &\times{{{\ddot{\varphi }}}_{1}} \left.\right.+\left( -\cos \varphi \sin {{\theta }_{1}}{{w}_{1}}+{{x}_{1}}\cos \varphi \cos {{\theta }_{1}} \right)\ddot{\theta }+\cos \varphi \cos {{\theta }_{1}}\frac{{{\partial }^{2}}{{w}_{1}}}{\partial {{t}^{2}}}+\left(\right. 2\sin {{\theta }_{1}}\sin \varphi {{w}_{1}}-2\\
  &\times{{x}_{1}}\cos{{\theta }_{1}}\sin\varphi){{{\dot{\theta }}}_{1}}\dot{\varphi }-2\sin \varphi \cos {{\theta }_{1}}\dot{\varphi }\frac{\partial {{w}_{1}}}{\partial t}+(-a\sin \varphi -\cos {{\theta }_{1}}\cos \varphi {{w}_{1}} -{{x}_{1}}\cos \varphi\sin {{\theta }_{1}})\\
  &\times{{{\dot{\varphi }}}^{2}}-2\cos \varphi \sin {{\theta }_{1}}{{{\dot{\theta }}}_{1}}\frac{\partial {{w}_{1}}}{\partial t}+\left( -\cos {{\theta }_{1}}\cos \varphi {{w}_{1}}-{{x}_{1}}\cos \varphi \sin {{\theta }_{1}} \right){{{\dot{\theta }}}_{1}}^{2} \left.\right]d{{x}_{1}}-\int_{0}^{L}{\rho A}[ \ddot{X}-\\
  &(a\cos \varphi+\cos {{\theta }_{2}}\sin \varphi {{w}_{2}}+{{x}_{2}}\sin {{\theta }_{2}}\sin \varphi)\ddot{\varphi } -\left( \cos \varphi \sin {{\theta }_{2}}{{w}_{2}}-{{x}_{2}}\cos \varphi \cos {{\theta }_{2}} \right){{{\ddot{\theta }}}_{2}}+\cos \varphi\\
  &\times\cos {{\theta }_{2}}\frac{{{\partial }^{2}}{{w}_{2}}}{\partial {{t}^{2}}} +\left(\right. 2\sin {{\theta }_{2}}\sin \varphi {{w}_{2}}-2{{x}_{2}}\cos {{\theta }_{2}}\sin \varphi  \left.\right){{{\dot{\theta }}}_{2}}\dot{\varphi }-2\sin \varphi \cos {{\theta }_{2}}\dot{\varphi }\frac{\partial {{w}_{2}}}{\partial t}+( a\sin \varphi\\
  &-\cos {{\theta }_{2}}\cos \varphi {{w}_{2}}-{{x}_{2}}\cos \varphi \sin {{\theta }_{2}}){{{\dot{\varphi }}}^{2}}-2\cos \varphi \sin {{\theta }_{2}}{{{\dot{\theta }}}_{2}}\frac{\partial {{w}_{2}}}{\partial t}-\left(\right. \cos {{\theta }_{2}}\cos \varphi {{w}_{2}}+{{x}_{2}}\cos \varphi\\
  &\times\sin {{\theta }_{2}} \left.\right){{{\dot{\theta }}}_{2}}^{2} \left.\right]d{{x}_{2}}=0 
\end{split}
\end{align}

\begin{align}
\label{eq25}
\begin{split}
   &\left( \frac{{{\tau }_{1}}+{{\tau }_{2}}}{{{r}_{w}}} \right)\sin \varphi +\cos \varphi {{F}_{s}}-\frac{2{{I}_{{{y}_{wheel}}}}}{{{r}_{w}}^{2}}\left( \sin 2\varphi \dot{Y}\dot{\varphi }+\cos 2\varphi \dot{X}\dot{\varphi } \right)-2\left(\right. \frac{{{I}_{{{y}_{wheel}}}}}{{{r}_{w}}^{2}}{{\sin }^{2}}\varphi +{{m}_{w}}\\
   &+{{m}_{base}}\left.\right)\ddot{Y}-\frac{{{I}_{{{y}_{wheel}}}}}{{{r}_{w}}^{2}}\ddot{X}\sin 2\varphi-\int_{0}^{L}{\rho A}\left[\right. \ddot{Y}+( a\sin \varphi +\cos {{\theta }_{1}}\cos \varphi {{w}_{1}}+{{x}_{1}}\cos \varphi \sin {{\theta }_{1}}) \\
   &\times\ddot{\varphi }+\left(\right. {{x}_{1}}\cos {{\theta }_{1}}\sin \varphi -\sin {{\theta }_{1}}\sin \varphi {{w}_{1}} \left.\right){{{\ddot{\theta }}}_{1}}+\cos {{\theta }_{1}}\sin \varphi \frac{{{\partial }^{2}}{{w}_{1}}}{\partial {{t}^{2}}}+\left(\right. 2{{x}_{1}}\cos {{\theta }_{1}}\cos \varphi -2\times\\
   &\cos \varphi \sin {{\theta }_{1}}{{w}_{1}})\dot{\varphi }{{{\dot{\theta }}}_{1}}+2\cos {{\theta }_{1}}\cos \varphi \dot{\varphi }\frac{\partial {{w}_{1}}}{\partial t}+\left(\right. a\cos \varphi -\cos {{\theta }_{1}}\sin \varphi {{w}_{1}}-{{x}_{1}}\sin {{\theta }_{1}}\sin \varphi  \left.\right){{{\dot{\varphi }}}^{2}}\\
  &-2\sin {{\theta }_{1}}\sin \varphi {{{\dot{\theta }}}_{1}}\frac{\partial {{w}_{1}}}{\partial t} +( -{{x}_{1}}\sin {{\theta }_{1}}\sin \varphi -\cos {{\theta }_{1}}\sin \varphi {{w}_{1}}){{{\dot{\theta }}}_{1}}^{2} ]d{{x}_{1}}-\int_{0}^{L}{}\rho A[ \ddot{Y}-(a\sin \varphi \\
  &-\cos {{\theta }_{2}}\cos \varphi {{w}_{2}}-{{x}_{2}}\cos \varphi\sin {{\theta }_{2}} \left.\right)\ddot{\varphi }{}+\left( {{x}_{2}}\cos {{\theta }_{2}}\sin \varphi -\sin {{\theta }_{2}}\sin \varphi {{w}_{2}} \right){{{\ddot{\theta }}}_{2}}+\cos {{\theta }_{2}}\sin \varphi \\
  &\frac{{{\partial }^{2}}{{w}_{2}}}{\partial {{t}^{2}}}+2\cos {{\theta }_{2}}\cos \varphi \dot{\varphi }\frac{\partial {{w}_{2}}}{\partial t}+\left(\right. 2{{x}_{2}}\cos {{\theta }_{2}}\cos \varphi-2\cos \varphi\sin {{\theta }_{2}}{{w}_{2}} \left.\right)\dot{\varphi }{{{\dot{\theta }}}_{2}}-2\sin {{\theta }_{2}}\sin \varphi {{{\dot{\theta }}}_{2}} \\
  &\times\frac{\partial {{w}_{2}}}{\partial t}+\left( -a\cos \varphi -{{w}_{2}}\cos {{\theta }_{2}}\sin \varphi -{{x}_{2}}\sin {{\theta }_{2}}\sin \varphi  \right){{{\dot{\varphi }}}^{2}}+\left(\right. -{{x}_{2}}\sin {{\theta }_{2}}\sin \varphi-\cos {{\theta }_{2}}\times \\
  &\sin \varphi {{w}_{2}} \left.\right){{{\dot{\theta }}}_{2}}^{2} \left.\right]d{{x}_{2}}=0  
\end{split}\\
\label{eq26}
\begin{split}
   & 2a\left( \frac{{{\tau }_{1}}-{{\tau }_{2}}}{{{r}_{w}}} \right)-(\frac{8{{a}^{2}}{{I}_{{{y}_{wheel}}}}}{{{r}_{w}}^{2}}+8{{a}^{2}}{{m}_{w}}+2{{a}^{2}}{{m}_{base}}+2{{I}_{{{z}_{wheel}}}}+2{{I}_{{{z}_{base}}}})\ddot{\varphi } +\frac{{{I}_{{{y}_{wheel}}}}}{{{r}_{w}}^{2}}\sin 2\varphi\\
   &\times{{{\dot{Y}}}^{2}}-\frac{{{I}_{{{y}_{wheel}}}}}{{{r}_{w}}^{2}}\sin 2\varphi {{{\dot{X}}}^{2}} +2\frac{{{I}_{{{y}_{wheel}}}}}{{{r}_{w}}^{2}}\dot{Y}\dot{X}\cos 2\varphi-\int_{0}^{L}{\rho A}[\left(\right. a\cos \varphi -\sin \varphi \cos {{\theta }_{1}}{{w}_{1}}-{{x}_{1}}\\
   &\sin \varphi \sin {{\theta }_{1}})\ddot{X} +a\cos {{\theta }_{1}}\frac{{{\partial }^{2}}{{w}_{1}}}{\partial {{t}^{2}}} +\left(\right. a\sin \varphi +{{x}_{1}}\cos \varphi \sin {{\theta }_{1}} +{{w}_{1}}\cos \varphi\cos {{\theta }_{1}} \left.\right)\ddot{Y} +\left(\right. {{a}^{2}}+\\
   &{{x}_{1}}^{2}+{{\cos }^{2}}{{\theta }_{1}}{{w}_{1}}^{2}-{{x}_{1}}^{2}{{\cos }^{2}}{{\theta }_{1}}+{{x}_{1}}\sin 2{{\theta }_{1}}{{w}_{1}})\ddot{\varphi }+(a{{x}_{1}}\cos {{\theta }_{1}}-a\sin {{\theta }_{1}}{{w}_{1}}){{{\ddot{\theta }}}_{1}} +\left(\right. 2{{\cos }^{2}}{{\theta }_{1}}\\
   & \times{{w}_{1}} +{{x}_{1}}\sin 2{{\theta }_{1}} )\dot{\varphi }\frac{\partial {{w}_{1}}}{\partial t}-2a\sin {{\theta }_{1}}{{{\dot{\theta }}}_{1}}\frac{\partial {{w}_{1}}}{\partial t}+(-a\cos {{\theta }_{1}}{{w}_{1}}-a{{x}_{1}}\sin {{\theta }_{1}} ){{{\dot{\theta }}}_{1}}^{2} +({{x}_{1}}^{2}\sin 2{{\theta }_{1}}\\
   &-\sin 2{{\theta }_{1}}{{w}_{1}}^{2}+2{{x}_{1}}\cos 2{{\theta }_{1}}{{w}_{1}} \left.\right)\dot{\varphi }{{{\dot{\theta }}}_{1}}\left.\right]d{{x}_{1}}-\int_{0}^{L}{\rho A}[-( a\cos \varphi +\sin \varphi \cos {{\theta }_{2}}{{w}_{2}}+{{x}_{2}}\sin \varphi\\
   &\times\sin {{\theta }_{2}})\ddot{X} +( {{a}^{2}}+{{x}_{2}}^{2}+{{\cos }^{2}}{{\theta }_{2}}{{w}_{2}}^{2}-{{x}_{2}}^{2}{{\cos }^{2}}{{\theta }_{2}}+{{x}_{2}}\sin 2{{\theta }_{2}}{{w}_{2}})\ddot{\varphi }-(a\sin \varphi -{{x}_{2}}\cos \varphi\\
   &\times\sin {{\theta }_{2}}-\cos \varphi \cos {{\theta }_{2}}{{w}_{2}})\ddot{Y}+\left(\right. -a{{x}_{2}}\cos {{\theta }_{2}} +a\sin {{\theta }_{2}}{{w}_{2}}){{{\ddot{\theta }}}_{2}}-a\cos {{\theta }_{2}}\frac{{{\partial }^{2}}{{w}_{2}}}{\partial {{t}^{2}}}+( a\cos {{\theta }_{2}}\\
   &\times{{w}_{2}}+a{{x}_{2}}\sin {{\theta }_{2}}){{{\dot{\theta }}}_{2}}^{2}+\left( 2{{\cos }^{2}}{{\theta }_{2}}{{w}_{2}}+{{x}_{2}}\sin 2{{\theta }_{2}} \right)\dot{\varphi }\frac{\partial {{w}_{2}}}{\partial t} +2a\sin {{\theta }_{2}}{{{\dot{\theta }}}_{2}}\frac{\partial {{w}_{2}}}{\partial t} +({{x}_{2}}^{2}\sin 2{{\theta }_{2}}\\
   &-\sin 2{{\theta }_{2}}{{w}_{2}}^{2}+2{{x}_{2}}{{w}_{2}}\cos 2{{\theta }_{2}} \left.\right)\dot{\varphi }{{{\dot{\theta }}}_{2}}\left.\right]d{{x}_{2}}=0
\end{split} \\
\label{eq27}
\begin{split}
 &-{{I}_{{{y}_{base}}}}{{{\ddot{\theta }}}_{1}}-\int_{0}^{L}{\rho A}\left[ \left( {{w}_{1}}^{2}+{{x}_{1}}^{2} \right){{{\ddot{\theta }}}_{1}}+\left( {{x}_{1}}\cos {{\theta }_{1}}\cos \varphi -\cos \varphi \sin {{\theta }_{1}}{{w}_{1}} \right)\ddot{X} \right.+{{x}_{1}}\frac{{{\partial }^{2}}{{w}_{1}}}{\partial {{t}^{2}}}+\\
 &(a{{x}_{1}}\cos {{\theta }_{1}}-a\sin {{\theta }_{1}}{{w}_{1}})\ddot{\varphi } +({{x}_{1}}\cos {{\theta }_{1}}\sin \varphi -\sin {{\theta }_{1}}\sin \varphi {{w}_{1}})\ddot{Y}+2{{w}_{1}}{{{\dot{\theta }}}_{1}}\frac{\partial {{w}_{1}}}{\partial t}\left(\right.\sin {{\theta }_{1}}{{w}_{1}}^{2}\\
 &\times\cos {{\theta }_{1}}-{{x}_{1}}^{2}\sin {{\theta }_{1}}\cos {{\theta }_{1}}-{{x}_{1}}\cos 2{{\theta }_{1}}{{w}_{1}} ){{{\dot{\varphi }}}^{2}}-g{{w}_{1}}\cos {{\theta }_{1}}-{{x}_{1}}\sin {{\theta }_{1}} ]d{{x}_{1}}=0
\end{split} \\
\label{eq28}
\begin{split}
  &-{{I}_{{{y}_{base}}}}{{{\ddot{\theta }}}_{2}}-\int_{0}^{L}{\rho A}[ \left( {{w}_{2}}^{2}+{{x}_{2}}^{2} \right){{{\ddot{\theta }}}_{2}}+\left( {{x}_{2}}\cos {{\theta }_{2}}\cos \varphi -\cos \varphi \sin {{\theta }_{2}}{{w}_{2}} \right)\ddot{X} +{{x}_{2}}\frac{{{\partial }^{2}}{{w}_{2}}}{\partial {{t}^{2}}} -\\
  &\left( a{{x}_{2}}\cos {{\theta }_{2}}-a\sin {{\theta }_{2}}{{w}_{2}} \right)\ddot{\varphi }+({{x}_{2}} \cos {{\theta }_{2}}\sin \varphi -\sin {{\theta }_{2}}\sin \varphi {{w}_{2}})\ddot{Y}+2{{w}_{2}}{{{\dot{\theta }}}_{2}}\frac{\partial {{w}_{2}}}{\partial t} +(\sin {{\theta }_{2}}\\
  &\times\cos {{\theta }_{2}}{{w}_{2}}^{2}-{{x}_{2}}^{2}\sin {{\theta }_{2}}\cos {{\theta }_{2}}-{{x}_{2}}\cos 2{{\theta }_{2}}{{w}_{2}}){{{\dot{\varphi }}}^{2}}-g{{w}_{2}}\cos {{\theta }_{2}}-g{{x}_{2}}\sin {{\theta }_{2}}]d{{x}_{2}}=0 
\end{split} \\
\label{eq29}
\begin{split}
&\rho A \left(\right. \frac{{{\partial }^{2}}{{w}_{1}}}{\partial {{t}^{2}}}+{{x}_{1}}{{{\ddot{\theta }}}_{1}}+a\cos {{\theta }_{1}}\ddot{\varphi }+\cos {{\theta }_{1}}\cos \varphi \ddot{X}+\cos {{\theta }_{1}}\sin \varphi \ddot{Y} -{{w}_{1}}{{{\dot{\theta }}}_{1}}^{2}+( -{{\cos }^{2}}{{\theta }_{1}}{{w}_{1}}-\\
&{{x}_{1}}\sin {{\theta }_{1}}\cos {{\theta }_{1}} ){{{\dot{\varphi }}}^{2}})+{{C}_{2}}\frac{\partial {{w}_{1}}}{\partial t}+EI\frac{{{\partial }^{4}}{{w}_{1}}}{\partial {{x}_{1}}^{4}}+\rho gA\sin {{\theta }_{1}}+\frac{{{\partial }^{2}}}{\partial {{x}_{1}}^{2}}\left( \frac{EIS(x)z{{d}_{31}}}{{{t}_{p}}}{{V}_{1}}(t) \right)=0 
\end{split} 
\end{align}

\begin{align}
\label{eq30}
\begin{split}
  &\rho A(\frac{{{\partial }^{2}}{{w}_{2}}}{\partial {{t}^{2}}}+{{x}_{2}}{{{\ddot{\theta }}}_{2}}-a\cos {{\theta }_{2}}\ddot{\varphi }+\cos {{\theta }_{2}}\cos \varphi \ddot{X}+\cos {{\theta }_{2}}\sin \varphi \ddot{Y} -{{w}_{2}}{{{\dot{\theta }}}_{2}}^{2}+\left(\right. -{{\cos }^{2}}{{\theta }_{2}}{{w}_{2}}-\\
  &{{x}_{2}}\sin {{\theta }_{2}}\cos {{\theta }_{2}} ){{{\dot{\varphi }}}^{2}}) +{{C}_{1}}\frac{\partial {{w}_{2}}}{\partial t}+EI\frac{{{\partial }^{4}}{{w}_{2}}}{\partial {{x}_{2}}^{4}}+\rho gA\sin {{\theta }_{2}}+\frac{{{\partial }^{2}}}{\partial {{x}_{2}}^{2}}\left( \frac{EIS(x)z{{d}_{31}}}{{{t}_{p}}}{{V}_{2}}(t) \right)=0  
\end{split}    
\end{align}

\end{strip}
where
\begin{equation}\label{eq31}
S(x)=H(x)-H(x-{{L}_{p}})
\end{equation}
and $H(x)$ is the Heaviside function. The boundary conditions of the above equations are represented as: 

\begin{equation}
 \label{eq32}
\left( \frac{{{\partial }^{2}}{{w}_{1}}}{\partial {{x}_{1}}^{2}}\delta \frac{\partial {{w}_{1}}}{\partial {{x}_{1}}}-\frac{{{\partial }^{3}}{{w}_{1}}}{\partial {{x}_{1}}^{3}}\delta {{w}_{1}} \right)_{0}^{L}=0   
\end{equation}

\begin{equation}
\label{eq33}
\left( \frac{{{\partial }^{2}}{{w}_{2}}}{\partial {{x}_{2}}^{2}}\delta \frac{\partial {{w}_{2}}}{\partial {{x}_{2}}}-\frac{{{\partial }^{3}}{{w}_{2}}}{\partial {{x}_{2}}^{3}}\delta {{w}_{2}} \right)_{0}^{L}=0
\end{equation}

To show the accuracy of Eqs. \eqref{eq24}-\eqref{eq33}, it is sufficient to omit the $\varphi$ and one of these beams. Then, the final equations will convert to\textcolor{blue}{\cite{mehrvarz2018modeling}}.

\section{Numerical simulation}\label{numericalsimulation}
As described in the previous section, the derived equations of motion are complex and have many nonlinear and coupled terms. In this section, a solution technique for these nonlinear equations is briefly presented. The natural frequencies of the beams need to be obtained first in order to solve the obtained equations of motion and extract the natural frequencies of the whole system. In order to extract the natural frequencies of the beams, the following undamped, unforced equations of motion and boundary conditions of the beams are to be used: 
\begin{equation}
 \label{eq34}
\rho A\frac{{{\partial }^{2}}w}{\partial {{t}^{2}}}+EI\frac{{{\partial }^{4}}w}{\partial {{x}^{4}}}=0   
\end{equation}

\begin{equation}
\label{eq35}
\left\{ \begin{matrix}
   w\left( 0,t \right)=0 \\ 
  \frac{\partial w}{\partial x}\left( 0,t \right)=0 \\ 
  \frac{{{\partial }^{2}}w}{\partial {{x}^{2}}}\left( L,t \right)=0 \\
  \frac{{{\partial }^{3}}w}{\partial {{x}^{3}}}\left( L,t \right)=0 \\
\end{matrix} \right.
\end{equation}

Here, the beams are assumed to have harmonic motions with frequency $\omega$ as\textcolor{blue}{\cite{jalili2009piezoelectric}}:

\begin{equation}\label{eq36}
w=W(x){{e}^{i\omega t}}
\end{equation}

Substituting Eq. \eqref{eq36} to \eqref{eq34}, results the below equation:

\begin{equation}\label{eq37}
-{{\rho }_{b}}{{A}_{b}}{{\omega }^{2}}W+E{{I}_{b}}{{W}''}''=0
\end{equation}

The general solution for \eqref{eq37} and its boundary conditions is given as

\begin{equation}\label{eq38}
W(x)={{a}_{1}}\cos \beta x+{{a}_{2}}\sin \beta x+{{a}_{3}}\cosh \beta x+{{a}_{4}}\sinh \beta x
\end{equation}
where
\begin{equation}\label{eq39}
{{\beta }^{2}}=\sqrt{\frac{\rho A}{EI}}\omega 
\end{equation}

Substituting boundary conditions \eqref{eq35} into eigenfunction \eqref{eq38}, results in the following set of equations:
\begin{equation}\label{eq40}
{{A}_{2\times 2}}{{X}_{2\times 1}}=0
\end{equation}
where
\begin{equation}\label{eq41}
{{X}_{2\times 1}}=\left[ \begin{matrix}
   {{a}_{1}}  \\
   {{a}_{2}}  \\
\end{matrix} \right]
\end{equation}

To obtain a nontrivial solution, the determinant of matrix $A$ in Eq. \eqref{eq40} should equal to zero. This equation gives the natural frequencies of the beams. The next step is to solve the equations of motion of the system using the assumed mode model expansion technique and the obtained natural frequencies. In this technique, the lateral displacements $w_1$ and $w_2$ are assumed as follows\textcolor{blue}{\cite{rao2007vibration}}:
\begin{equation}\label{eq42}
\begin{matrix}
  & {{w}_{1}}=\sum\limits_{i=1}^{\infty }{{{W}_{i}}({{x}_{1}}){{q}_{1i}}(t)} \\ 
 & {{w}_{2}}=\sum\limits_{i=1}^{\infty }{{{W}_{i}}({{x}_{2}}){{q}_{2i}}(t)} \\ 
\end{matrix}
\end{equation}
where $q_{1i} (t)$ and $q_{2i} (t)$ are the generalized coordinates for the bending of the beams and $W_i (x_1)$ and $W_i (x_2)$ are the mode shapes of a fixed-free beam. These functions are defined as: 
\begin{equation}\label{eq43}
\begin{split}
  &{{W}_{i}}({{x}_{i}})={{{A}'}_{1}}(sin({{\beta }_{n}}{{x}_{i}})-sinh({{\beta }_{n}}{{x}_{i}}))\\
&+{{{A}'}_{2}}(cos({{\beta }_{n}}{{x}_{i}})-cosh({{\beta }_{n}}{{x}_{i}}))  
\end{split}
\end{equation}
where $A'_1$ and $A'_2$ are two tunable and dependent coefficients as 
\begin{equation}\label{eq44}
{{{A}'}_{2}}=-\frac{(sin({{\beta }_{n}}L)-sinh({{\beta }_{n}}L))}{(cos({{\beta }_{n}}L)+cosh({{\beta }_{n}}L))}{{{A}'}_{1}}
\end{equation}
with $\beta_n$ being defined as follows for each mode:
\begin{equation}\label{eq45}
{{\beta }_{n}}^{4}=\frac{{{\rho }_{b}}{{\omega }_{1n}}^{2}}{E{{I}_{b}}}
\end{equation}

The equations of motion can be obtained by substituting Eq. \eqref{eq42} to Eqs. \eqref{eq24}-\eqref{eq30} and also multiplying Eq. \eqref{eq43} to Eqs. \eqref{eq29} and \eqref{eq30}, and then integrating the obtained equations over $0$ to $L$. The final equations represent $2n+5$ DOF of the system. 

\section{Simulation Results} \label{simulationresults}
\begin{figure}
\centering
\includegraphics[width=80mm]{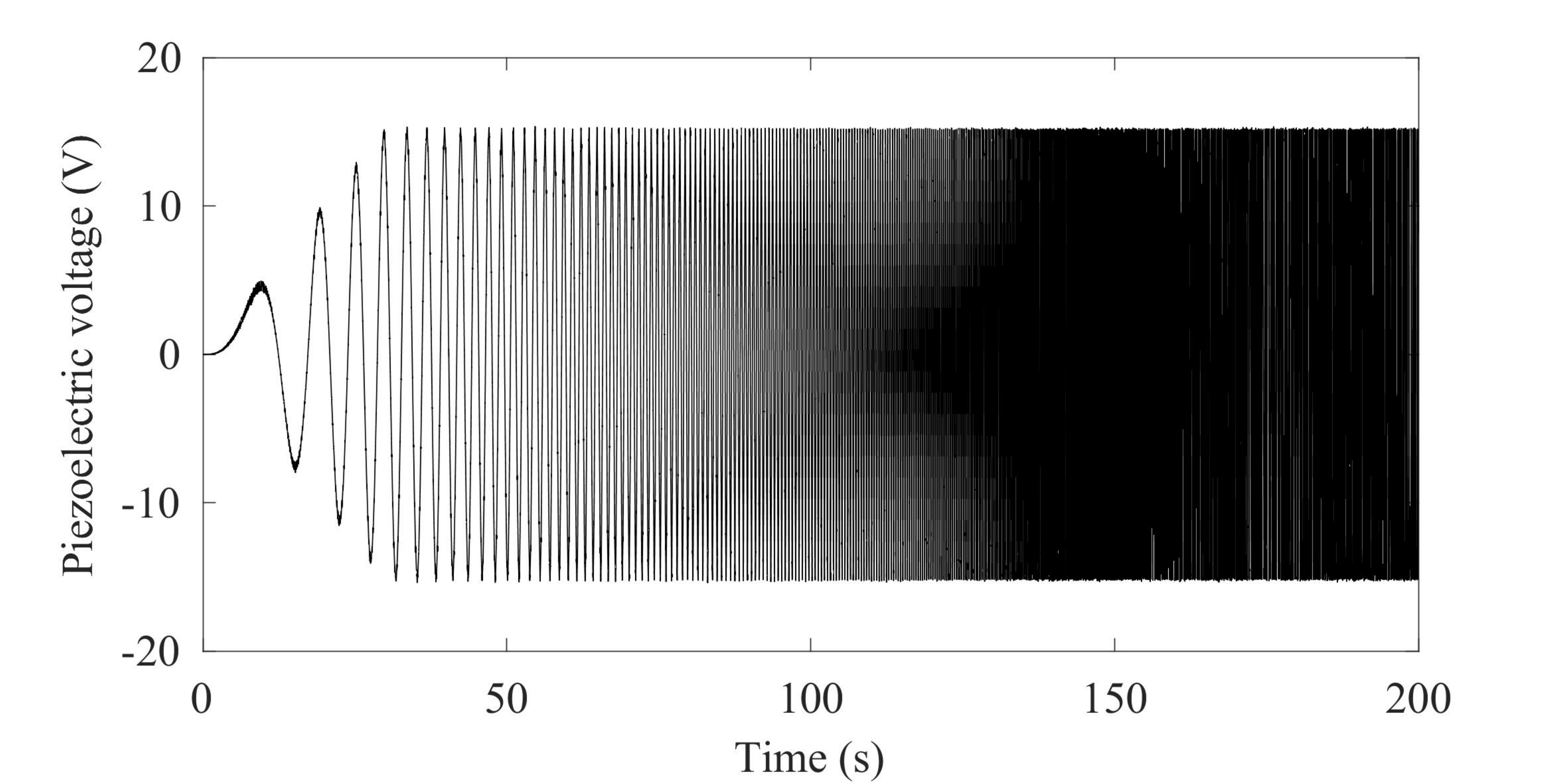}
\caption{The swipe frequency input voltage.}
\label{fig4} 
\end{figure}

\begin{figure}
\centering
\includegraphics[width=80mm]{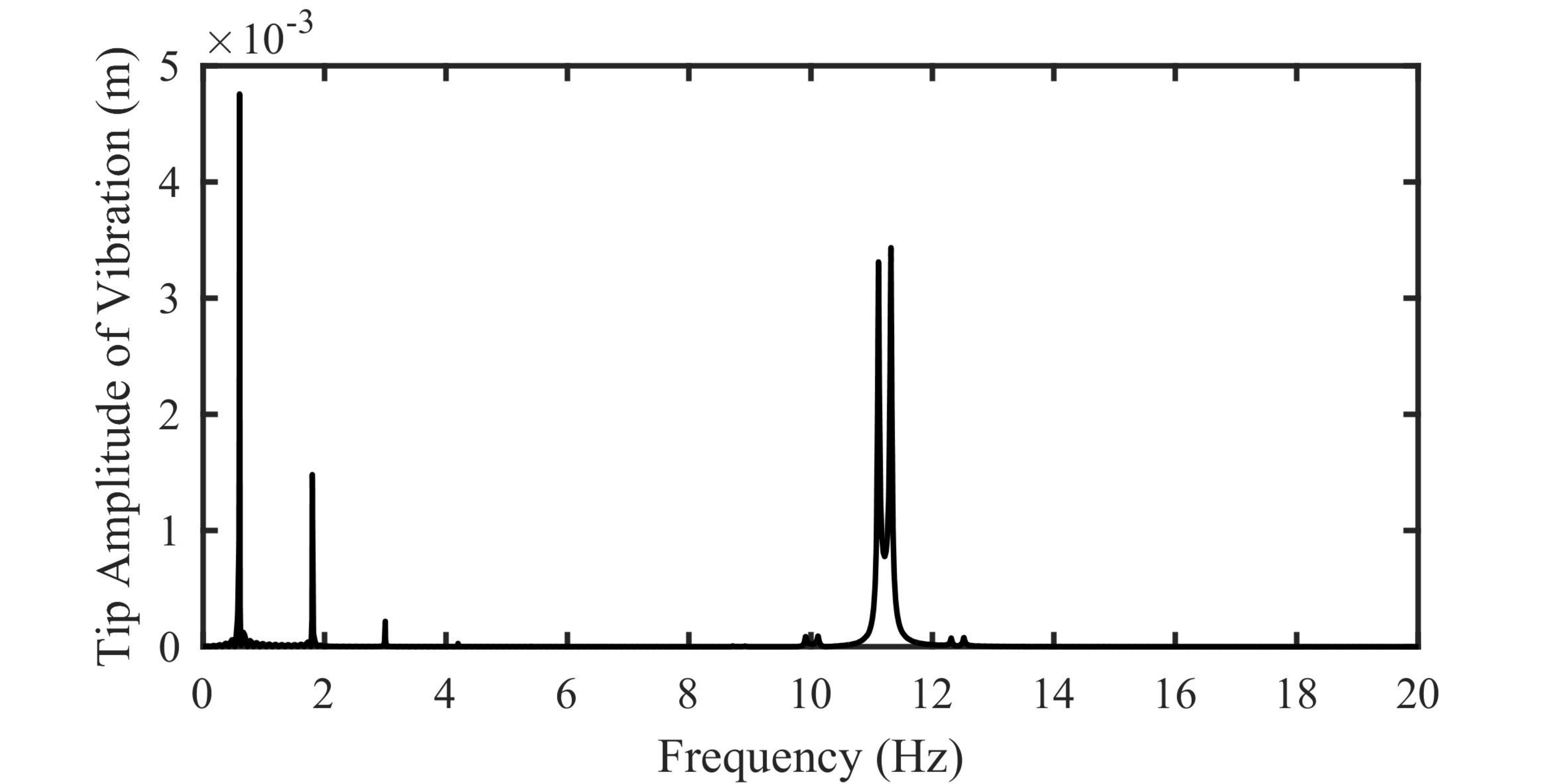}
\caption{The FFT analysis of the system.}
\label{fig5}
\end{figure}

To investigate the dynamic behavior of the system, the equations of motions are solved numerically in Matlab and the results are presented in different scenarios. The numerical values of the physical parameters are presented in \textcolor{blue}{{Table }\ref{table1}} and only two modes are considered for both beams. First, a sweep frequency input, shown in \textcolor{blue}{{Fig.}\ref{fig4}}, is applied to the system to extract the natural frequencies of the system. This standard input provides a fairly uniform spectral excitation and gives the modes of the systems\textcolor{blue}{\cite{banazadeh2017identification}}. \textcolor{blue}{{Fig.}\ref{fig5}} shows the spectral analysis of time history of the system and its natural frequencies.

In the next scenario, the lateral movement of the system in $X$-direction is presented by applying two same unit step input voltages to the piezoelectric actuators. It is assumed that the robot just moves forward. Hence, the below relation is set between the beams' rotation angle and the DC motors' voltage.
\begin{equation}\label{eq46}
\begin{matrix}
   \left\{ \begin{matrix}
   {{V}_{ai}}={{10}^{-2}}{{\theta }_{i}} & \theta >0  \\
   {{V}_{ai}}=0 & \theta <0  \\
\end{matrix} \right. & ,i=1,2  \\
\end{matrix}
\end{equation}

\begin{figure}
\centering
\includegraphics[width=80mm]{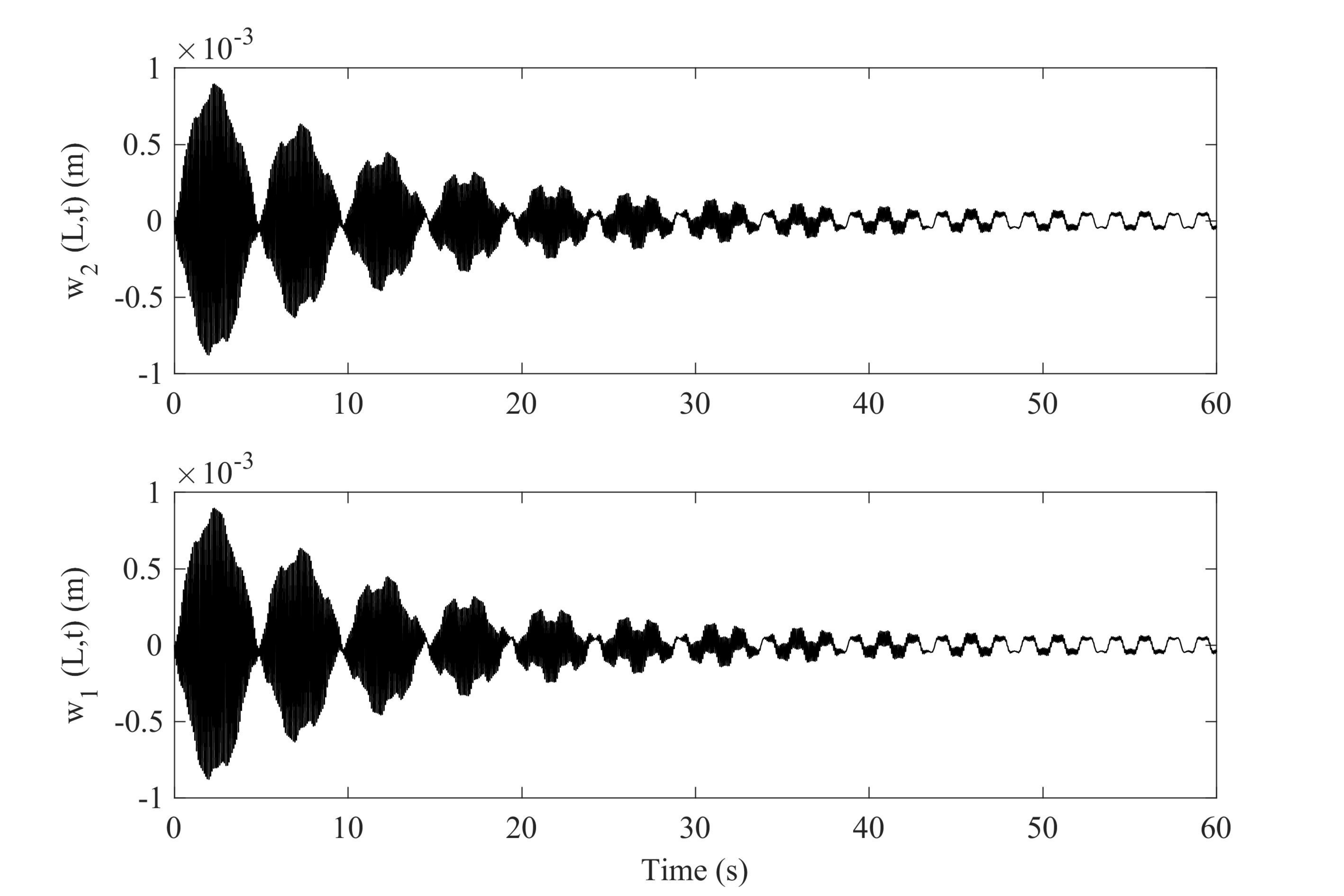}
\caption{Tip deflection of the beams $w_1 (L,t)$ and $w_2 (L,t)$ with unit step input.}
\label{fig6} 
\end{figure}

\begin{figure}
\centering
\includegraphics[width=80mm,scale=0.5]{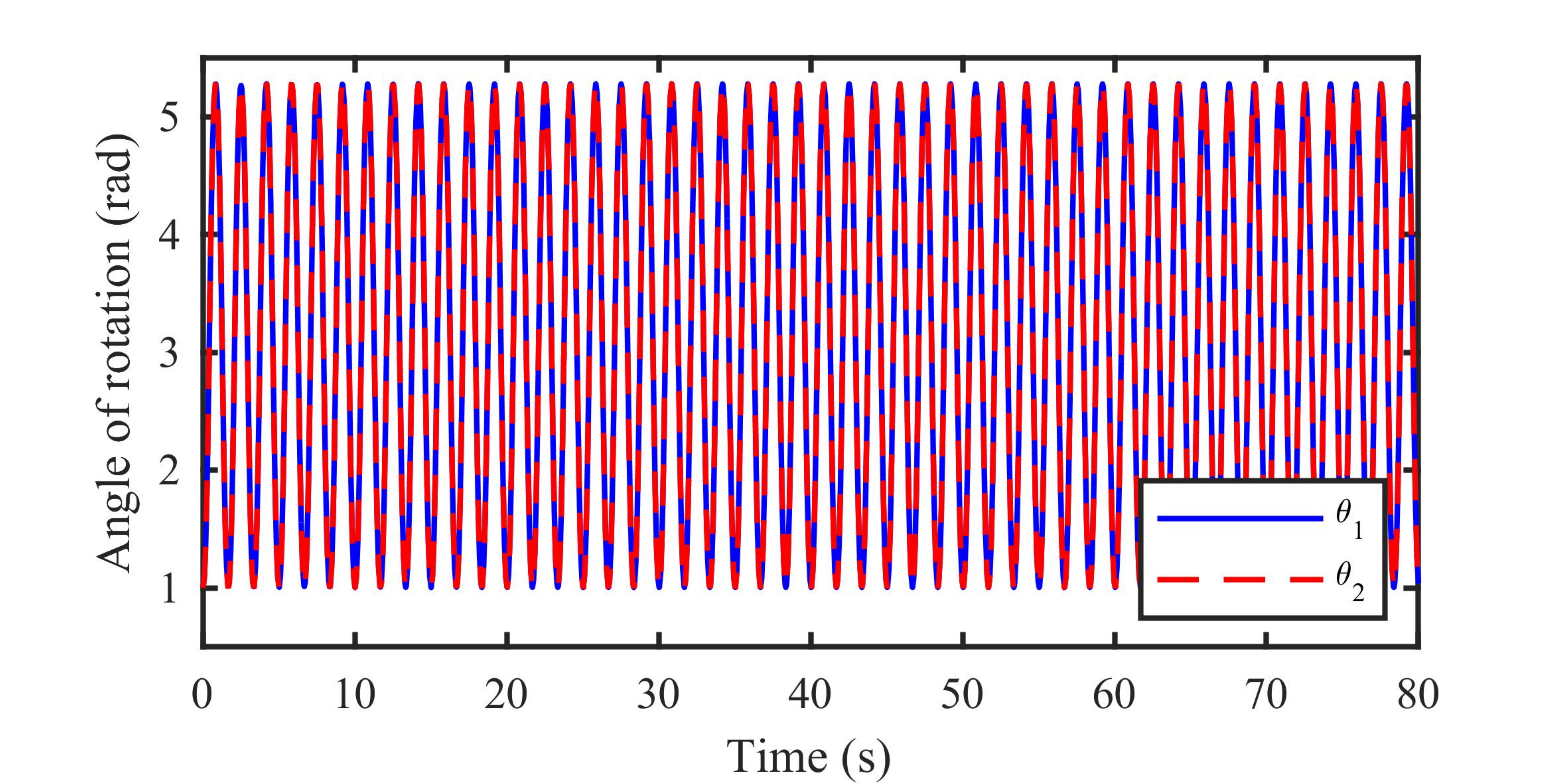}
\caption{Angular rotations of the beams $\theta_1$ and $\theta_2$  around the robot's base with unit step input.}
\label{fig7}
\end{figure}

\begin{figure}
\centering
\includegraphics[width=80mm]{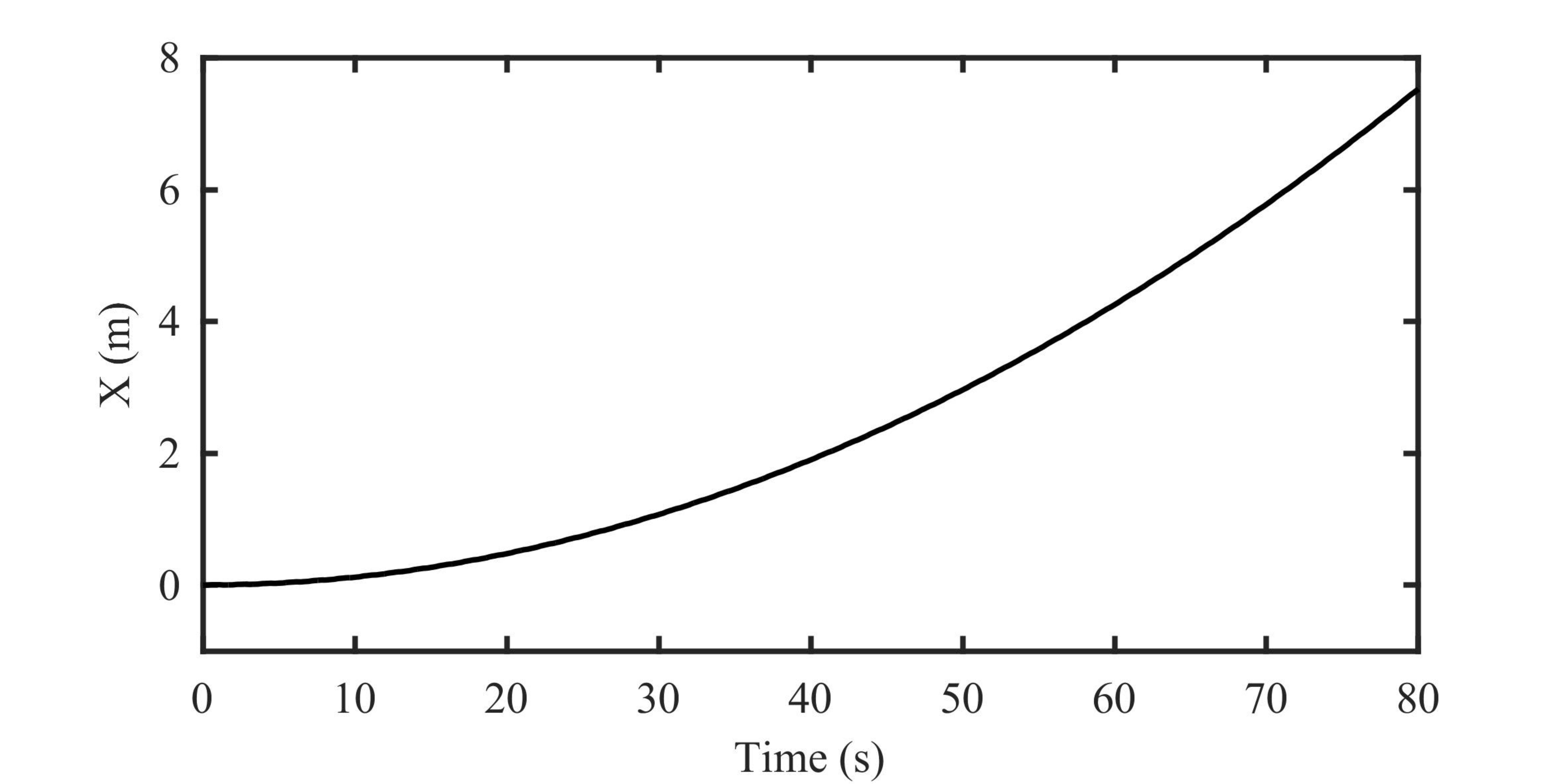}
\caption{Displacement of the robot in $X$-direction with step input.}
\label{fig8}
\end{figure}

\textcolor{blue}{{Fig.}\ref{fig6}} shows the lateral vibrations of the beams. These lateral vibrations, as shown in \textcolor{blue}{{Fig.}\ref{fig7}}, make the beam rotation $\theta$ variations and consequently produce the DC motors voltages. As considered in the previous section, these voltages cause the external torques and the lateral movement in $X$-direction. The displacement of the robot is shown in \textcolor{blue}{{Fig.}\ref{fig8}}. Because of nonlinearity of the system, $\varphi$ and $Y$ are not zero but their size are small with respect to both $X$ and $\theta_i$ and considered negligible.

\begin{figure}
\centering
\includegraphics[width=80mm]{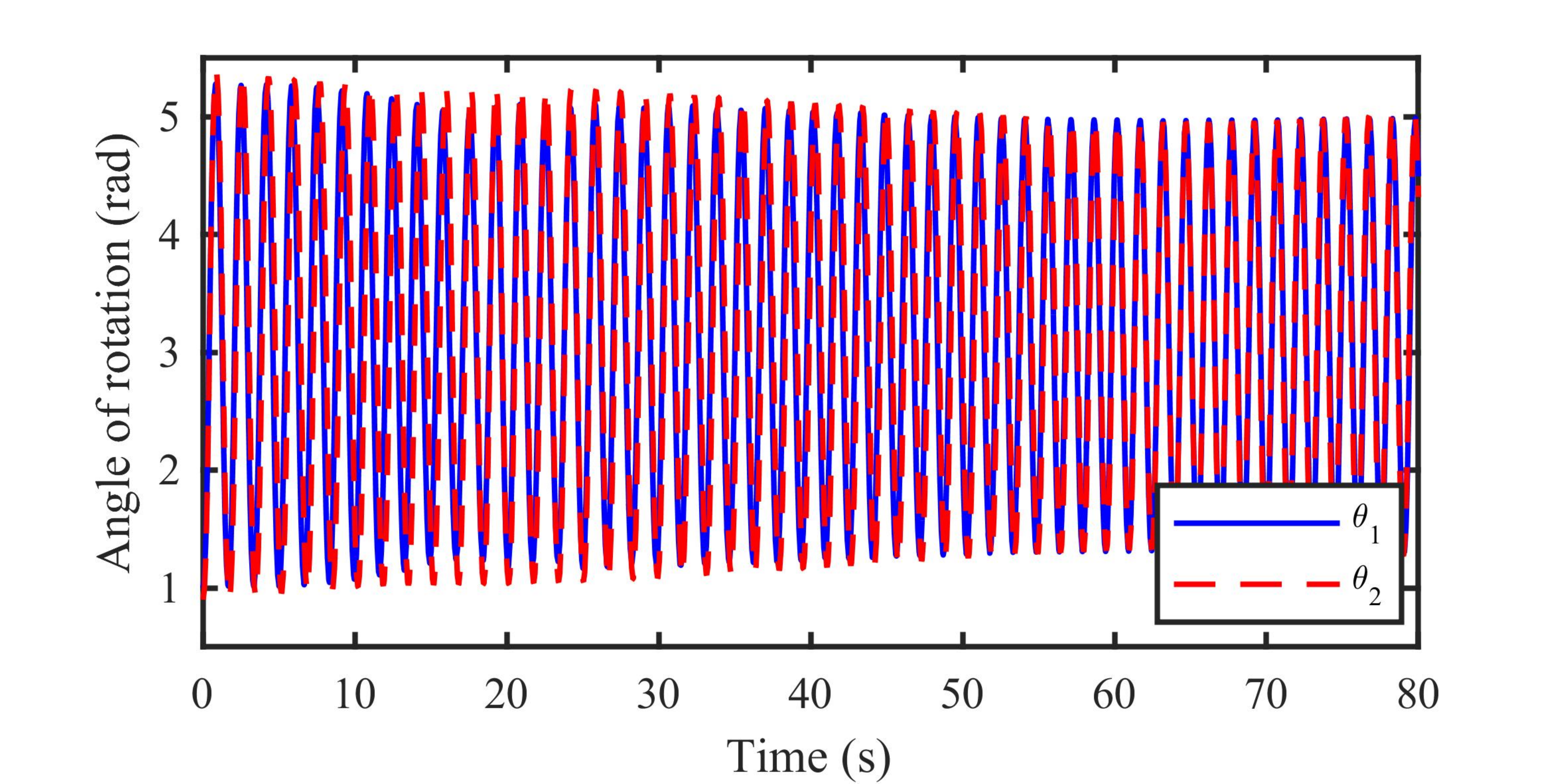}
\caption{Angular rotations of the beams $\theta_1$ and $\theta_2$ around the robot's base with $v_1=H(t)$ and $v_1=\frac{9}{10} H(t)$.}
\label{fig9} 
\end{figure}

\begin{figure}
\centering
\includegraphics[width=80mm]{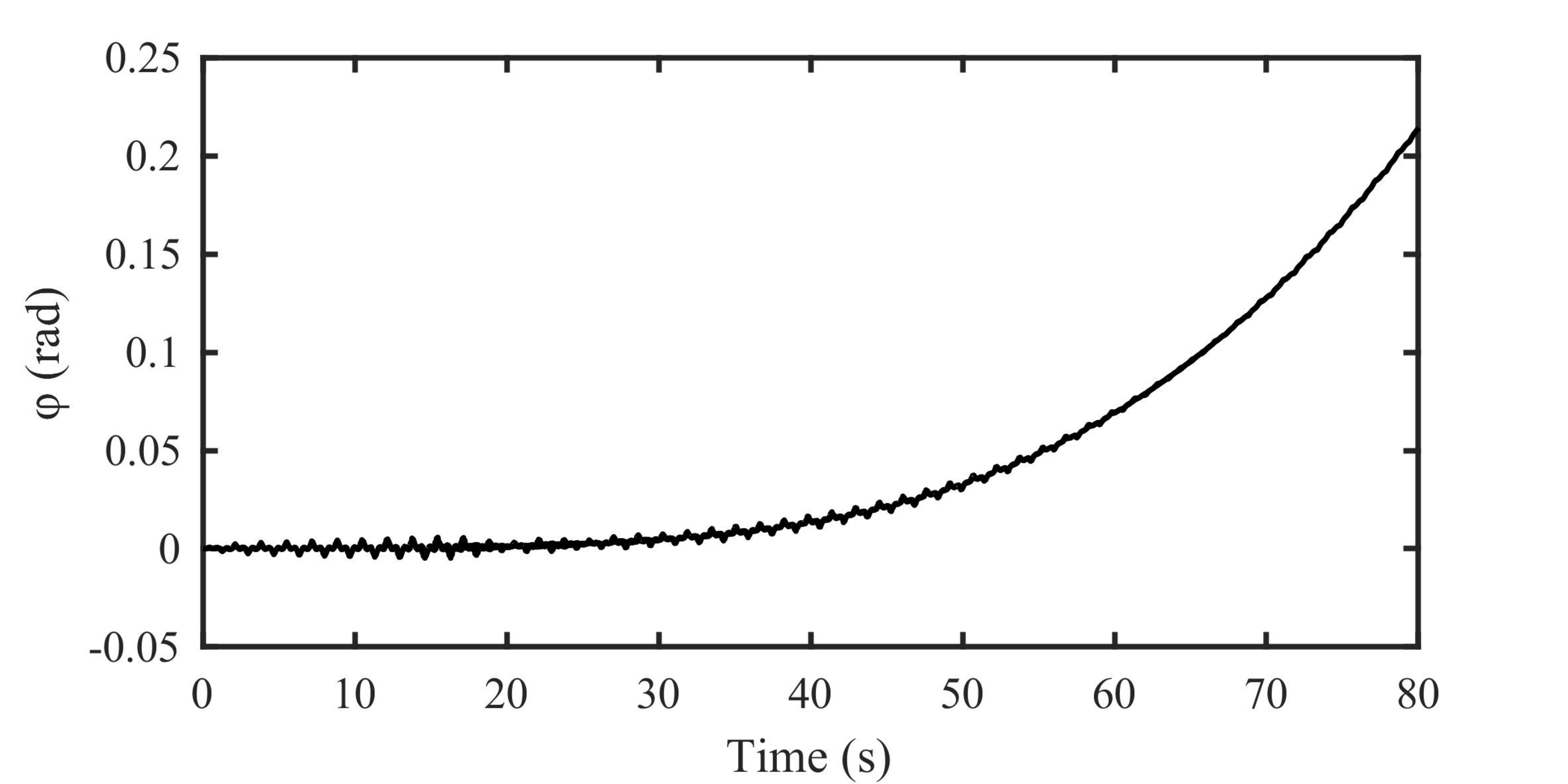}
\caption{Angular rotation $\varphi$ of the robot around the $z$-axis.}
\label{fig10}
\end{figure}

\begin{figure}
\centering
\includegraphics[width=80mm]{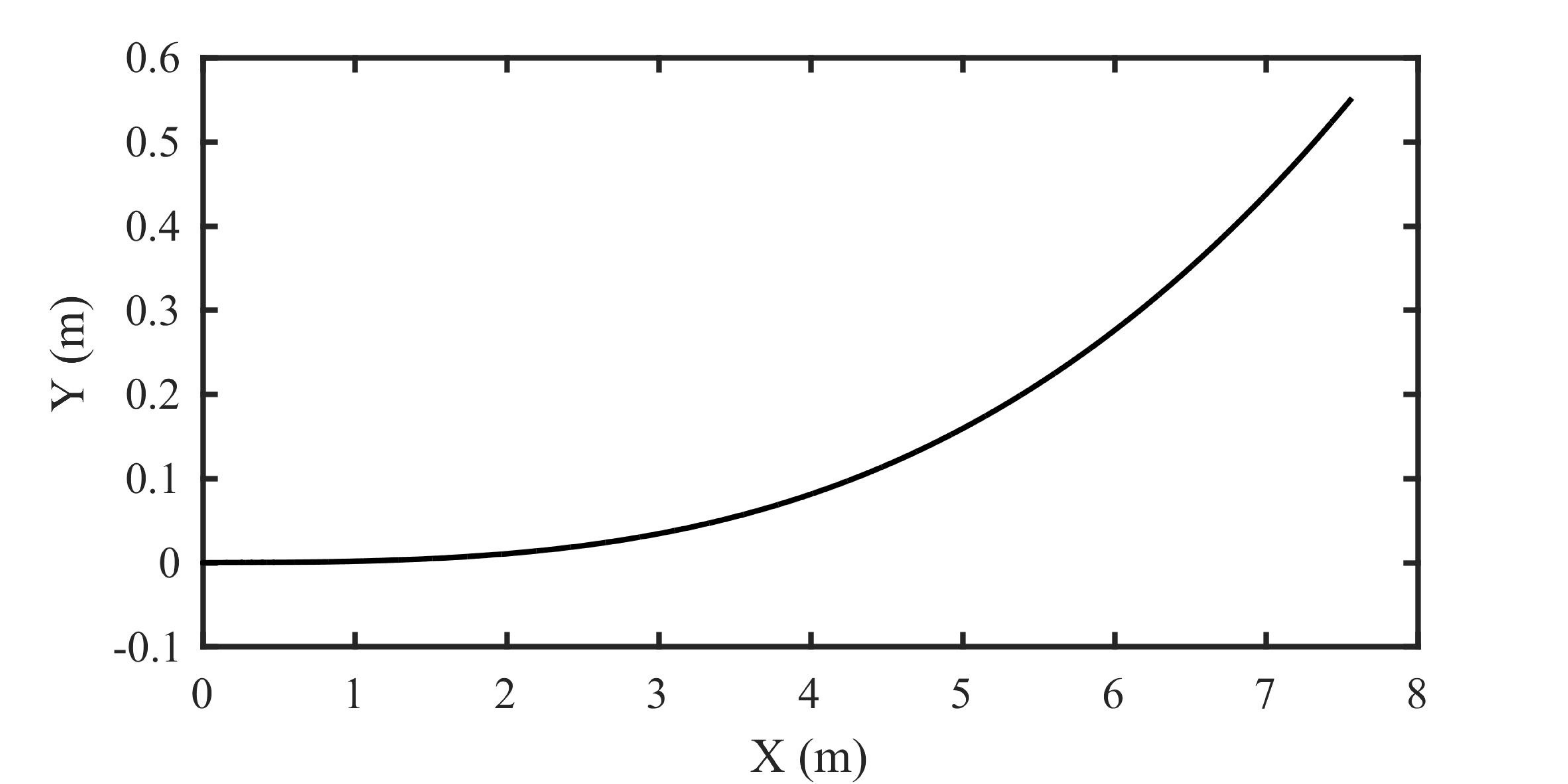}
\caption{Displacement $X$ and $Y$ of the robot in $XY$-plane.}
\label{fig11}
\end{figure}

The final scenario is rotation test in the $X$-$Y$ plane. In this test, two step input voltages with different amplitudes are applied to the piezoelectric actuators. This makes different beams' lateral vibrations and bases' rotational angles, and consequently different torques applied to the wheels. This causes different angular velocities in the wheels and produces rotational angle $\varphi$ in the robot system. \textcolor{blue}{{Fig.}\ref{fig9}} shows the rotation angles $\theta_1$ and $\theta_2$. The rotation and the path of the robot are also shown in \textcolor{blue}{{Fig.}\ref{fig10}} and \textcolor{blue}{{Fig.}\ref{fig11}}.

\begin{table*}
\centering
\caption{The system parameters.}
{\footnotesize
\begin{tabular}
{llll}
\hline \hline
Parameter & Value & Parameter & Value\\ \hline \hline
Beam length (mm) & 271.46 &
Piezo layer shear modulus (GPa) & 5.515 \\
Beam thickness (mm) & 0.5 &
Piezo layer elastic modulus (GPa) & 30.33 \\
Beam width (mm) & 25.65 &
Piezo layer density (Kg/m3) & 5440 \\
Beam elastic modulus (GPa) & 70 &
First flexural damping ratio (\%) & 0.0058 \\
Beam shear modulus (GPa) & 30 &
Second flexural damping ratio (\%) & 0.015 \\
Piezo layer length (mm) & 38 &
PZT layer width (mm) & 23 \\
Piezo layer thickness (mm) & 0.3 &
Beam density (Kg/m3) & 2700 \\ \hline
\end{tabular}
}
\label{table1}
\end{table*}

\section{Conclusions}\label{conclusion}
A new configuration for the conventional two-wheeled inverted pendulum system was presented in this paper. The developed model had $2n+5$ DOF and its main purpose was to simulate two cantilever beams and piezoelectric actuators on a moving base as a new robot which can move in-plane. The governing equations of motion were obtained by employing the extended Hamilton's Principle. The derivation steps of these equations were presented in detail. The obtained model indicates that these systems have several coupled and nonlinear terms in their dynamics. To investigate the dynamic behavior of the system, these complex equations were solved numerically and the natural frequencies of the system were extracted. Finally, the result of two different tests including the lateral and circular movements were presented. 

\bibliographystyle{IEEEtran}
\bibliography{main.bib}
\flushend

\end{document}